\documentclass[12pt]{article}%
\usepackage{amssymb}
\usepackage{amsfonts}
\usepackage{amsmath}
\usepackage{graphicx}
\usepackage{amscd}%
\providecommand{\U}[1]{\protect\rule{.1in}{.1in}}
\newtheorem{theorem}{Theorem}[section]
\newtheorem{proposition}[theorem]{Proposition}
\newtheorem{lemma}[theorem]{Lemma}
\newtheorem{corollary}[theorem]{Corollary}

\begin{document}

\title{A new reconstruction method in integral geometry}
\author{V. P. Palamodov\\Tel Aviv University}
\date{}
\maketitle

\textbf{Abstract: }A general method for analytic inversion of geometric
integral transforms is proposed.



\section{Introduction}

There are numerous applications of integral geometry to various problem of
image reconstruction, in particular to X-ray, emission, magnetic resonance,
wave, acoustic, thermo and photoacoustic, scattering, Doppler tomographies as
well as to radar technique and texture analysis. We propose here a general
method of analytic inversion formulas for various transformations in integral
geometry. Several examples for the case of integral transforms in a plane
domain are given.

\section{Curves and integrals}

Let $X$ and $\Sigma$ be smooth 2-manifolds and $\Phi$ be a smooth real
function defined in $X\times\Sigma$ such that $\mathrm{d}_{x}\Phi\neq0.$ For
any point $\sigma\in\Sigma$ the set $F\left(  \sigma\right)  =\{x\in
X;\Phi\left(  x,\sigma\right)  =0\}$ is a smooth curve in $X.$ We call $\Phi$
\textit{generating} function of the curve family $\{F\left(  \sigma\right)
,\sigma\in\Sigma\}.$

Let $\mathrm{d}S$ be the Riemannian area form; we define a Funk-Radon
transform generated by $\Phi$ by%
\[
M_{\Phi}f\left(  \sigma\right)  =\int\delta\left(  \Phi\left(  x;\sigma
\right)  \right)  f\mathrm{d}S=\int_{F\left(  \sigma\right)  }\frac
{f\mathrm{d}S}{\mathrm{d}_{x}\Phi\left(  x;\sigma\right)  }%
\]
for continuous functions $f$ compactly supported\ in $X.$ The quotient
$f\mathrm{d}S/\mathrm{d}_{x}\Phi$ denotes an arbitrary 1-form $q$ such that
$\mathrm{d}_{x}\Phi\wedge q=f\mathrm{d}S.$ It is defined up to a term
$h\mathrm{d}_{x}\Phi$ where $h$ is a continuous function. An orientation of a
curve $F\left(  \sigma\right)  $ is defined by means of the form
$\mathrm{d}_{x}\Phi$ and the integral of $q$ over $F\left(  \sigma\right)  $
is uniquely defined. Let $\mathrm{g}$ be a Riemannian metric in a manifold
$X;$ the form $\mathrm{d}_{\mathrm{g}}s=\sqrt{g\left(  \mathrm{d}x\right)  }$
is the Riemannian line element.$\ $We have then%
\[
M_{\Phi}f\left(  \sigma\right)  =\int_{F\left(  \sigma\right)  }%
\frac{f\mathrm{d}_{\mathrm{g}}s}{\left\vert \nabla_{\mathrm{g}}\Phi\left(
x;\sigma\right)  \right\vert }%
\]
where $\left\vert \nabla_{\mathrm{g}}a\right\vert \doteq\sqrt{\mathrm{g}%
\left(  \mathrm{d}a\right)  }$ is the Riemannian gradient of a function $a$.
Suppose that the gradient factorizes through $X$ and $\Sigma$ that is
$\left\vert \nabla_{\mathrm{g}}\Phi\left(  x,\sigma\right)  \right\vert
=m\left(  x\right)  \mu\left(  \sigma\right)  $ for some positive continuous
functions $m$ in $X$ and $\mu$ in $\Sigma.\ $It follows that data of the
Funk-Radon transform is equivalent to data of Riemannian curve integrals%
\[
R_{\mathrm{g}}f\left(  \sigma\right)  =\int_{F\left(  \sigma\right)
}f\mathrm{d}_{\mathrm{g}}s,\sigma\in\Sigma
\]
since $Rf\left(  \sigma\right)  =\mu\left(  \sigma\right)  M_{\Phi}\left(
mf\right)  \left(  \sigma\right)  .$ The reconstruction problem of a function
$f$ from Riemannian integrals $R_{\mathrm{g}}f$ is then reduced to inversion
of the operator $M_{\Phi}f.$

\section{Main theorem}

We assume now that $\Sigma=\mathbb{R}\times S^{1}$ and a generating function
is linear in the first argument: $\Phi\left(  x;\lambda,\varphi\right)
=\lambda+\psi\left(  x,\varphi\right)  ,$ $\lambda\in\mathbb{R},\varphi\in
S^{1}=\{0\leq\varphi<2\pi\}.\ $We call this function \textit{regular
}if\ (\textbf{i}) $\mathrm{d}_{x}\psi\wedge\mathrm{d}_{x}\psi_{\varphi
}^{\prime}\neq0$ in $X\times S^{1}$ and (\textbf{ii}) there are no conjugated
points, that is the equations $\psi\left(  x,\varphi\right)  =\psi\left(
y,\varphi\right)  $ and $\psi_{\varphi}^{\prime}\left(  x,\varphi\right)
=\psi_{\varphi}^{\prime}\left(  x,\varphi\right)  $ are fulfilled for no
$x\neq y\in X,$ $\varphi\in S^{1}$.

\begin{theorem}
\label{R}Let $\Phi\left(  x;\lambda,\varphi\right)  =\lambda+\psi\left(
x,\varphi\right)  $ be a regular generating function in $X\times\Sigma$ such
that the function $\psi$ is analytic in $X$ and the principal value integral%
\begin{equation}
N\left(  x,y\right)  \doteq\left(  P\right)  \int_{0}^{2\pi}\frac
{\mathrm{d}\varphi}{\left(  \psi\left(  x,\varphi\right)  -\psi\left(
y,\varphi\right)  \right)  ^{2}}\doteq\operatorname{Re}\int\frac
{\mathrm{d}\varphi}{\left(  \psi\left(  x,\varphi\right)  -\psi\left(
y,\varphi\right)  \pm i0\right)  ^{2}} \label{9}%
\end{equation}
vanishes for any $x\neq y\in X.$ Then the reconstruction formula%
\begin{align}
f\left(  x\right)   &  =-\frac{1}{4\pi^{2}D\left(  x\right)  }\left(
P\right)  \int_{0}^{2\pi}\int_{\mathbb{R}}\frac{M_{\Phi}f\left(
\lambda,\varphi\right)  }{\Phi^{2}\left(  x;\lambda,\varphi\right)
}\mathrm{d}\lambda\mathrm{d}\varphi\label{3}\\
D\left(  x\right)   &  =\frac{1}{2\pi}\int_{0}^{2\pi}\frac{\mathrm{d}\varphi
}{\left\vert \nabla_{\mathrm{g}}\psi\left(  x,\varphi\right)  \right\vert
^{2}} \label{33}%
\end{align}
holds for an arbitrary function $f\in L_{2}\left(  X\right)  _{\mathrm{comp}}%
$. The integral (\ref{3}) converges in $L_{2}\left(  X\right)  _{\mathrm{loc}%
}$.
\end{theorem}

\textbf{Remark 1.} A more invariant form of (\ref{3}) is a reconstruction of
the density$t.$%
\[
f\left(  x\right)  \mathrm{d}S=-\frac{1}{4\pi^{2}}\frac{\mathrm{d}S}{D\left(
x\right)  }\left(  P\right)  \int_{0}^{2\pi}\int_{\mathbb{R}}\frac{M_{\Phi
}f\left(  \lambda,\varphi\right)  }{\Phi^{2}\left(  x;\lambda,\varphi\right)
}\mathrm{d}\lambda\mathrm{d}\varphi
\]
Here the quotient $\mathrm{d}S/D\left(  x\right)  $ depends only of the
conformal class of the Riemannian metric $\mathrm{g}.$

\textbf{Remark 2. }G. Beylkin studied "the generalized Radon transform" in
$\mathbb{R}^{n}\mathbb{\ }$\cite{B} which coincides with our Funk-Radon
transform $M_{\Phi}$ in a Euclidean space (see also the last Section). He
constructed a Fourier integral operator which serves as parametrix for this
transform.\ Beylkin's parametrix gives a high frequency approximation to a
solution but not an exact inversion due to limitation of the method of FIO.
The explicit form of the nucleus $N\left(  x,y\right)  $ does not appear in
Beylkin's approach but is important in our method. Note that if the nucleus
does not vanish it is still a smooth function in the complement to the
diagonal, since condition (\textbf{ii) }keeps the number of real zeros of
$\psi\left(  x,\varphi\right)  -\psi\left(  y,\varphi\right)  $ constant. Then
the right-hand side of (\ref{3}) is a parametrix.

\textbf{Proof of Theorem}.\textbf{\ }To simplify our arguments we assume that
$X$ is an open set in $\mathbb{R}^{2}$ and $\mathrm{g}$ is the Euclidean
metric, the area form is denoted $\mathrm{d}x$. Let $q$ be a 1-form in $X$
such that $\mathrm{d}\psi\wedge q=\mathrm{d}x;$ we have $Mf\left(
\lambda,\varphi\right)  =\int_{\lambda+\psi=0}fq.$ For an arbitrary $x\in X$
and any function $f$ that vanishes in a neighborhood of $x $ we find%
\begin{align*}
\int_{\Sigma}\frac{M_{\Phi}f\left(  \lambda,\varphi\right)  \mathrm{d}%
\lambda\mathrm{d}\varphi}{\Phi^{2}\left(  x;\lambda,\varphi\right)  }  &
=\int_{S^{1}}\mathrm{d}\varphi\int_{\mathbb{R}}\left(  \int_{\lambda
+\psi\left(  y,\varphi\right)  =0}f\left(  y\right)  q\left(  y\right)
\right)  \frac{\mathrm{d}\lambda}{\left(  \lambda+\psi\left(  x,\varphi
\right)  \right)  ^{2}}\\
&  =-\int_{X}\left(  \int_{0}^{2\pi}\frac{\mathrm{d}\varphi}{\left(
\psi\left(  x,\varphi\right)  -\psi\left(  y,\varphi\right)  \right)  ^{2}%
}\right)  f\left(  y\right)  q\left(  y\right)  \wedge\mathrm{d}\psi\left(
y,\varphi\right) \\
&  =\int_{X}N\left(  x,y\right)  f\left(  y\right)  \mathrm{d}y
\end{align*}
since $\mathrm{d}\lambda=-\mathrm{d}\psi$ as $\lambda+\psi\left(
y,\varphi\right)  =0$ and $\mathrm{d}\lambda\wedge q=-\mathrm{d}y$ by
definition of $q$. It follows that the function $N$ is the off-diagonal
nucleus of $I$ and it vanishes since of the assumption. Therefore the nucleus
$I\left(  x,y\right)  $ of the operator $I$ (see (\ref{16}) is supported in
the diagonal and by Lemma \ref{CO} below we have $I\left(  x,y\right)
=a\left(  x\right)  \delta_{x}\left(  y\right)  $ for some continuous function
$a.$ To calculate this function at a point $x_{0}\in X$ we apply $I$ to a
density $h_{\varepsilon}~$that is equal to $\mathrm{d}x$ in $\varepsilon
$-neighborhood $X_{\varepsilon}$ of $x_{0}$ and $h_{\varepsilon}=0$ in the
complement:
\begin{equation}
Ih_{\varepsilon}\left(  x_{0}\right)  =\int_{X_{\varepsilon}}\mathrm{d}%
x\int_{0}^{2\pi}\frac{\mathrm{d}\varphi}{\left(  \psi\left(  x,\varphi\right)
-\psi\left(  x_{0},\varphi\right)  \right)  ^{2}}=\int_{0}^{2\pi}%
\mathrm{d}\varphi\int_{X_{\varepsilon}}\frac{\mathrm{d}x}{\left(  \psi\left(
x,\varphi\right)  -\psi\left(  x_{0},\varphi\right)  \right)  ^{2}}\nonumber
\end{equation}

\begin{lemma}
\label{D}We have for arbitrary $\varphi\in S^{1},\ x_{0}\in X$ and small
$\varepsilon$%
\[
D\left(  x_{0},\varphi\right)  \doteq\int_{X_{\varepsilon}}\frac{\mathrm{d}%
x}{\left(  \psi\left(  x,\varphi\right)  -\psi\left(  x_{0},\varphi\right)
\right)  ^{2}}=-\frac{2\pi}{\left\vert \nabla\psi\left(  x_{0},\varphi\right)
\right\vert ^{2}}+O\left(  \varepsilon\right)
\]
where $O\left(  \varepsilon\right)  \leq C\varepsilon$ where $C$ does not
depend of $\varphi.$
\end{lemma}

Integrating over $\varphi$ yields the equation $D\left(  x_{0}\right)
=-\int\left\vert \nabla\psi\left(  x_{0},\varphi\right)  \right\vert
^{-2}\mathrm{d}\varphi+O\left(  \varepsilon\right)  .$ This implies (\ref{33})
since $D\left(  x_{0}\right)  $ does not depend of $\varepsilon.$

\textbf{Proof of Lemma.} Replacing $\psi$ by $\psi/\left\vert \nabla
\psi\left(  x_{0},\varphi\right)  \right\vert $ we can assume that $\left\vert
\nabla\psi\left(  x_{0},\varphi\right)  \right\vert =1.$ Shift the origin to
the point $x_{0}$ and suppose first that $\psi\left(  x,\varphi\right)  $ is
linear in $x.$ By an orthogonal transformation we can assume that $\psi\left(
x,\varphi\right)  =x_{1}.$ In this case we have%
\begin{align*}
D\left(  x_{0},\varphi\right)   &  =\mathrm{\operatorname{Re}}\int
_{X_{\varepsilon}}\frac{\mathrm{d}x}{\left(  x_{1}+i0\right)  ^{2}%
}=-\mathrm{\operatorname{Re}}\int_{-\varepsilon}^{\varepsilon}\mathrm{d}%
x_{2}\left.  \frac{1}{x_{1}+i0}\right\vert _{-\sqrt{\varepsilon^{2}-x_{2}^{2}%
}}^{\sqrt{\varepsilon^{2}-x_{2}^{2}}}\\
&  =-2\int_{-\varepsilon}^{\varepsilon}\frac{\mathrm{d}x_{2}}{\sqrt
{\varepsilon^{2}-x_{2}^{2}}}=-\frac{2\pi}{\left\vert \nabla\psi\left(
x_{0},\varphi\right)  \right\vert ^{2}}%
\end{align*}
In the general case we can write $\psi^{2}\left(  x,\varphi\right)  =x_{1}%
^{2}-\rho\left(  x\right)  ,$ where $\rho\left(  x\right)  =O\left(
\left\vert x\right\vert ^{3}\right)  $ is analytic in $3\varepsilon
$-neighborhood of the origin if $\varepsilon$ is sufficiently small. We have%
\[
\frac{1}{\psi^{2}}=\frac{1}{x_{1}^{2}}+\sigma,\ \sigma\doteq\sum_{n=1}%
^{\infty}\frac{\rho^{n}}{x_{1}^{2n+2}}%
\]
and
\begin{equation}
\int_{X_{\varepsilon}}\frac{\mathrm{d}x}{\left(  \psi+i0\right)  ^{2}}%
=\int_{X_{\varepsilon}}\frac{\mathrm{d}x}{\left(  x_{1}+i0\right)  ^{2}}%
+\int_{X_{\varepsilon}}\sigma\left(  x_{1}+i0,x_{2}\right)  \mathrm{d}x
\label{10}%
\end{equation}
The first integral in the right-hand side is equal to $2\pi$ and shall show
the second integral is small.\ The field%
\[
v\left(  x\right)  =\left(  \left\vert x_{2}\right\vert +\sqrt{\varepsilon
^{2}-\left\vert x\right\vert ^{2}},-\frac{x_{2}}{\left\vert x_{2}\right\vert
}x_{1}\right)
\]
is defined and continuous for $x\in X_{\varepsilon},x_{2}\neq0.$ For any
$t,\ 0\leq t\leq1$ we define 2-chains
\begin{align*}
V_{t,\pm}  &  =\left\{  z=z\left(  x,t\right)  \doteq\frac{x+itv\left(
x\right)  }{\sqrt{1-t^{2}}},\ \pm x_{2}\geq0,\ \left\vert x\right\vert
=\varepsilon\right\} \\
W_{t}  &  =\left\{  z=\frac{1}{\sqrt{1-t^{2}}}\left(  x_{1}+it\sqrt
{\varepsilon^{2}-x_{1}^{2}},istx_{1}\right)  ,\ \left\vert x_{1}\right\vert
\leq\varepsilon,\ -1\leq s\leq1\right\}
\end{align*}
oriented by the orientation of $X_{\varepsilon}.$ The boundary of the chain
$V_{t}=V_{t,+}+V_{t,-}+W_{t}$ is the image of $\partial X_{\varepsilon}$ under
the map $x\mapsto z\left(  x,t\right)  .$ The function $\sigma$ has a
holomorphic continuation to the domain $Z_{+}=\left\{  z\in\mathbb{C}%
^{2};\left\vert z\right\vert <3\varepsilon,y_{1}>0\right\}  $. The form
$\sigma\mathrm{d}z$ is holomorphic in $Z_{+}$ hence $\mathrm{d}\sigma
\wedge\mathrm{d}z=0$. Consider an open 3-chain $Z_{\delta}=\cup_{t=\delta
}^{1/2}V_{t}\subset Z$ where $0<\delta\leq1/2.$ It is contained in $Z_{+}$
since $y_{1}\doteq\operatorname{Im}z_{1}=t\left(  \left\vert x_{2}\right\vert
+\sqrt{\varepsilon^{2}-\left\vert x\right\vert ^{2}}\right)  >0$ and
$\left\vert v\right\vert \leq2\varepsilon.$ By Stokes' we have%
\[
\int_{\partial Z_{\delta}}\sigma\mathrm{d}z=\int_{Z_{\delta}}\mathrm{d}%
\sigma\wedge\mathrm{d}z=0
\]
The boundary of $Z_{\delta}$ is equal to $V_{1/2}-V_{\delta}+S+T$ where%
\begin{align*}
S  &  =\left\{  z=\frac{x+itv\left(  x\right)  }{\sqrt{1-t^{2}}},\ \left\vert
x\right\vert =\varepsilon,\ \delta\leq t\leq1/2\right\}  ,\ \\
T  &  =\left\{  z=\frac{\left(  \varepsilon,ist\varepsilon\right)  }%
{\sqrt{1-\left(  st\right)  ^{2}}},\varepsilon=\pm1,\delta\leq t\leq1/2,-1\leq
s\leq1\right\}
\end{align*}
Therefore%
\[
\int_{V_{1/2}}\sigma\mathrm{d}z-\int_{V_{\delta}}\sigma\mathrm{d}x+\int
_{S+T}\sigma\mathrm{d}z=0
\]
For any point $z\in S$ we have%
\[
z^{2}\doteq z_{1}^{2}+z_{2}^{2}=\left(  1-t^{2}\right)  ^{-1}\left(
\left\vert x\right\vert ^{2}+2it\left\langle x,v\right\rangle -t^{2}\left\vert
v\right\vert ^{2}\right)  =\left(  1-t^{2}\right)  ^{-1}\left(  \left\vert
x\right\vert ^{2}-t^{2}\left\vert x\right\vert ^{2}\right)  =\left\vert
x\right\vert ^{2}=\varepsilon^{2}%
\]
since $\varepsilon^{2}-\left\vert x\right\vert ^{2}=0$ and $\left\langle
x,v\right\rangle =0$ on $\partial X_{\varepsilon}.$ The same equation holds in
$T.$ It follows that the chain $S+T$ is contained in the complex algebraic
curve $z^{2}=\varepsilon^{2}$. This implies that the integral of the
holomorphic form $\sigma\mathrm{d}z$ over $S+T$ vanishes. This yields for an
arbitrary $\delta>0$%
\[
\int_{V_{\delta}}\sigma\mathrm{d}z=\int_{V_{1/2}}\sigma\mathrm{d}z=\sum
_{n=1}^{\infty}\int_{V_{1/2}}\frac{\rho^{n}}{z_{1}^{2n+2}}%
\]
The left-hand side tends to the second term of (\ref{10}) as $\delta
\rightarrow0.$ The right-hand side can be estimated from above by the sum%
\[
\sum_{n=1}^{\infty}\left\vert \int_{V_{1/2}}\frac{\rho^{n}\mathrm{d}z}%
{z_{1}^{2n+2}}\right\vert \leq\int_{V_{1/2}}\mathrm{d}S\sum M_{\varepsilon
}^{n}\max_{z\in V_{1/2}}\left\vert z_{1}\right\vert ^{-2n-2}%
\]
where
\[
M_{\varepsilon}=\max_{Z}\left\vert \rho\right\vert \leq C\varepsilon^{3}%
\]
and $\mathrm{d}S$ is the Euclidean area element. The area of the chain $V_{1}$
is estimated by$\ \mathrm{const~}\varepsilon^{2}$ since $\left\vert \nabla
v\right\vert $ is bounded in $X_{\varepsilon}.$ For we have any point $z\in
V_{1}$
\[
\left\vert z_{1}\right\vert ^{2}=x_{1}^{2}+y_{1}^{2}\geq\frac{4}{3}\left(
x_{1}^{2}+\frac{1}{4}\left(  x_{2}^{2}+\varepsilon^{2}-\left\vert x\right\vert
^{2}\right)  \right)  \geq\frac{1}{3}\varepsilon^{2}%
\]
This yields%
\[
\left\vert \int_{V_{1/2}}\sigma\mathrm{d}z\right\vert \leq C_{0}\sum
_{1}^{\infty}\frac{\varepsilon^{2}\left(  3M_{\varepsilon}\right)  ^{n}%
}{\varepsilon^{2\left(  n+1\right)  }}\leq C_{0}\sum_{1}^{\infty}%
C^{n}\varepsilon^{n}=C_{0}\frac{C\varepsilon}{1-C\varepsilon}=O\left(
\varepsilon\right)
\]
if $\varepsilon<C^{-2}.$ This completes proofs of Lemma \ref{D} and Theorem
\ref{R}. $\blacktriangleright$

\begin{lemma}
\label{CO}The integral transform
\begin{equation}
If\left(  x\right)  =\left(  P\right)  \int_{0}^{2\pi}\int_{\mathbb{R}}%
\frac{M_{\Phi}f\left(  \lambda,\varphi\right)  }{\Phi^{2}\left(
x;\lambda,\varphi\right)  }\mathrm{d}\lambda\mathrm{d}\varphi\label{16}%
\end{equation}
is a continuous operator $L_{2}\left(  X\right)  _{\mathrm{comp}}\rightarrow
L_{2}\left(  X\right)  _{\mathrm{loc}}.$
\end{lemma}

\textbf{Proof of Lemma.} The condition (\textbf{i}) can be written in the
form$\ J_{x;\xi}\left(  \Phi\right)  \neq0$ in $F$ where
\begin{equation}
J_{x;\lambda,\varphi}\left(  \Phi\right)  =\det\left(
\begin{array}
[c]{ccc}%
\frac{\partial^{2}\Phi}{\partial x_{1}\partial\lambda} & \frac{\partial
^{2}\Phi}{\partial x_{1}\partial\varphi} & \frac{\partial\Phi}{\partial x_{1}%
}\\
&  & \\
\frac{\partial^{2}\Phi}{\partial x_{2}\partial\lambda} & \frac{\partial
^{2}\Phi}{\partial x_{2}\partial\varphi} & \frac{\partial\Phi}{\partial x_{2}%
}\\
&  & \\
\frac{\partial\Phi}{\partial\lambda} & \frac{\partial\Phi}{\partial\varphi} &
0
\end{array}
\right)  \label{15}%
\end{equation}
$x_{1},x_{2}$ are coordinates in $X.$ According to \cite{H}, Theorem 25.3.1
the map $M_{\Phi}:L_{2}\left(  X\right)  _{\mathrm{comp}}\rightarrow
L_{2}\left(  \Sigma\right)  _{\mathrm{loc}}$ is a Fourier integral operator of
order $-1/2$ since the corresponding Lagrange manifold is locally a graph of a
canonical transformation (see details in \cite{P2}; the bundles $\Omega
^{1/2}\left(  X\right)  $ and $\Omega^{1/2}\left(  \Sigma\right)  $ are
trivial in our case). The image of this operator is contained in $L_{2}\left(
\Sigma\right)  _{\mathrm{comp}}$ since for an arbitrary function $f\in
L_{2}\left(  X\right)  _{\mathrm{comp}}$ the support of\textrm{\ }$M_{\Phi}f$
is contained in $\Theta=\pi_{\Sigma}\pi_{X}^{-1}\left(  \mathrm{supp~}%
f\right)  $ where $\pi_{X}$ and $\pi_{\Sigma}$ are projections of $F$ to $X$
and $\Sigma$ respectively. The projection $\pi_{X}$ is obviously
proper.\thinspace The integral transform
\[
Sh\left(  x\right)  =\int_{\Sigma}\frac{h\left(  \lambda,\varphi\right)
\mathrm{d}\lambda\mathrm{d}\varphi}{\Phi^{2}\left(  x;\lambda,\varphi\right)
}%
\]
can be written by means of a Fourier integral operator%
\[
Sh\left(  x\right)  =-\frac{1}{2}\int_{\Sigma}\int_{\mathbb{R}}\exp\left(
\imath t\Phi\left(  x,\lambda,\varphi\right)  \right)  \left\vert t\right\vert
h\left(  \lambda,\varphi\right)  \mathrm{d}t\mathrm{d}\lambda\mathrm{d}\varphi
\]
with the phase function $\Phi$. The corresponding Lagrange manifold coincides
with that of $M_{\Phi}$ with source and target spaces interchanged. The order
of $S$ equals $1/2$ due to the factor $\left\vert t\right\vert $. It follows
that the composition $I=SM_{\Phi}$ is well defined as a PDO of order
$0.\blacktriangleright$

\section{Integrals of rational trigonometric functions}

We study now the condition $N\left(  x,y\right)  =0$. A function%
\[
t\left(  \varphi\right)  =\sum_{m=0}^{k}a_{m}\cos m\varphi+b_{m}\sin m\varphi
\]
is called trigonometric polynomial of order $k$ if $\left\vert a_{k}%
\right\vert ^{2}+\left\vert b_{k}\right\vert ^{2}>0.$ Any trigonometric
polynomial is $2\pi$-periodic and is well-defined and holomorphic in the
cylinder $\mathbb{C}/2\pi\mathbb{Z}.$ It always has $2k$ zeros in the
cylinder. If a polynomial is real the number of real zeros is even.

\begin{lemma}
\label{T}Let $t\left(  \varphi\right)  ,s\left(  \varphi\right)  $ be real
trigonometric polynomials such that $\mathrm{\deg~}s<n~\mathrm{\deg~}t$ for a
natural $n$ and all zeros of $t$ are real and simple. Then%
\begin{equation}
\left(  P\right)  \int_{0}^{2\pi}\frac{s\left(  \varphi\right)  }{t^{n}\left(
\varphi\right)  }\mathrm{d}\varphi\doteq\frac{1}{2}\int_{0}^{2\pi}%
\frac{s\left(  \varphi\right)  }{\left(  t\left(  \varphi\right)  +i0\right)
^{n}}\mathrm{d}\varphi+\frac{1}{2}\int_{0}^{2\pi}\frac{s\left(  \varphi
\right)  }{\left(  t\left(  \varphi\right)  -i0\right)  ^{n}}\mathrm{d}%
\varphi=0 \label{7}%
\end{equation}

\end{lemma}

\textbf{Proof.} The function $r\left(  \zeta\right)  =s\left(  \zeta\right)
t^{-n}\left(  \zeta\right)  $ is meromorphic for $\zeta=\varphi+i\tau
\in\mathbb{C}/2\pi\mathbb{Z}$ \ and has no nonreal poles since of the
assumption. Fix a small positive $\varepsilon$ we choose a real continuous
function $\lambda=\tau\left(  \varphi\right)  $ defined on the circle
$0\leq\varphi\leq2\pi$ that vanishes except for $\varepsilon$-neighborhood of
the zero set of $t$ and
\[
\tau\left(  \varphi\right)  =\mathrm{sgn~}t^{\prime}\left(  \alpha\right)
\sqrt{\varepsilon^{2}-\left(  \varphi-\alpha\right)  ^{2}}%
\]
in $\varepsilon$-neighborhood of each zero $\alpha$ of $t.$ We have $t\left(
\alpha+i\tau\left(  \alpha\right)  0\right)  =t\left(  \alpha\right)  +i0$ for
any zero $\alpha$ and
\[
\int_{0}^{2\pi}\frac{s\left(  \varphi\right)  }{\left(  t\left(
\varphi\right)  +i0\right)  ^{n}}\mathrm{d}\varphi=\int_{0}^{2\pi}%
\frac{s\left(  \zeta\right)  }{t^{n}\left(  \zeta\right)  }\mathrm{d}%
\zeta,\ \zeta\left(  \varphi\right)  =\varphi+i\tau\left(  \varphi\right)
\]
This yields%
\[
\int_{0}^{2\pi}\frac{s\left(  \varphi\right)  }{\left(  t\left(
\varphi\right)  +i0\right)  ^{n}}\mathrm{d}\varphi+\int_{0}^{2\pi}%
\frac{s\left(  \varphi\right)  }{\left(  t\left(  \varphi\right)  -i0\right)
^{n}}\mathrm{d}\varphi=\int_{\Gamma\cup\bar{\Gamma}}\frac{s\left(
\zeta\right)  }{t^{n}\left(  \zeta\right)  }\mathrm{d}\zeta
\]
where $\Gamma=\left\{  \zeta=\zeta\left(  \varphi\right)  \right\}  .$ We can
write $\Gamma\cup\bar{\Gamma}=\Gamma_{+}\cup\Gamma_{-}$ where $\Gamma_{\pm
}=\left\{  \zeta=\varphi\pm i\left\vert \tau\left(  \varphi\right)
\right\vert \right\}  .$ Finally%
\[
\int_{\Gamma\cup\bar{\Gamma}}\frac{s\left(  \zeta\right)  }{t^{n}\left(
\zeta\right)  }\mathrm{d}\zeta=\int_{\Gamma_{+}}\frac{s\left(  \zeta\right)
}{t^{n}\left(  \zeta\right)  }\mathrm{d}\zeta+\int_{\Gamma_{-}}\frac{s\left(
\zeta\right)  }{t^{n}\left(  \zeta\right)  }\mathrm{d}\zeta=0
\]
since $t\left(  \varphi+i\tau\right)  \neq0$ for $\tau\neq0$ and\ the function
$st^{-n}$ tends to zero at infinity.\ $\blacktriangleright$

The following formula can be useful for calculation of $D\left(  x\right)  :$

\begin{lemma}
\label{I}If $t$ and $s$ are real trigonometric polynomials,$\ \mathrm{\deg
\ }s<\deg t=k$ and $t$ has no real zeros, then%
\[
\int_{0}^{2\pi}\frac{s\left(  \varphi\right)  }{t\left(  \varphi\right)
}\mathrm{d}\varphi=\operatorname{Re}\left(  2\pi i\sum_{t\left(  \varphi
_{m}\right)  =0,\operatorname{Im}\varphi_{m}>0}\frac{s\left(  \varphi
_{m}\right)  }{t^{\prime}\left(  \varphi_{m}\right)  }\right)
\]
where the sum is taken over $k$ zeros of $t$ with positive (or negative)
imaginary part.
\end{lemma}

For a proof we apply the Residue theorem for the form $s\left(  \varphi
\right)  \mathrm{d}\varphi/t\left(  \varphi\right)  $ in the upper
half-cylinder $\mathbb{C}_{+}/2\pi\mathbb{Z}.$

\section{Radon's and Funk's reconstructions}

We apply Theorem \ref{R} to recover some known and unknown formulas for
regular curve families.

\textbf{Radon's }formula\textbf{. }Take a generating function $\Phi\left(
x;\lambda,\omega\right)  =\lambda-\left\langle x,\mathbf{e}\left(
\varphi\right)  \right\rangle $ in $\mathbb{R}^{2}\times\Sigma$ where
$\mathbf{e}\left(  \varphi\right)  =\left(  \cos\varphi,\sin\varphi\right)  $.
The classical formula of Radon-John%
\begin{equation}
f\left(  x\right)  =-\frac{1}{2\pi^{2}}\int_{0}^{\pi}\int_{-\infty}^{\infty
}\frac{g\left(  \lambda,\varphi\right)  \mathrm{d}\lambda\mathrm{d}\varphi
}{\left(  \lambda-\left\langle x,\mathbf{e}\left(  \varphi\right)
\right\rangle \right)  ^{2}} \label{8}%
\end{equation}
coincides in this case with (\ref{3}) since integral data $g=R$ $f$ are equal
to $Mf\left(  \lambda,\varphi\right)  $ because of $\left\vert \nabla
\psi\right\vert =1.$ The equation $N\left(  x,y\right)  =0$ immediately
follows from Lemma \ref{T}. The coefficient (\ref{33}) is reduced
to$\ -1/2\pi^{2}$ due to symmetry $g\left(  -\lambda,\varphi+\pi\right)
=g\left(  \lambda,\varphi\right)  $.

\textbf{Funk's }formula provides reconstruction of an even function $f$ in the
unit sphere $S^{2}$ from integrals\ of $f$ over the family $\Sigma\ $of big
circles. Take the unit hemisphere $X=\{x\in E^{3},\left\vert x\right\vert
=1,x_{0}\geq0\}$ and a generating function $\Phi\left(  x;\lambda
,\varphi\right)  =\lambda+\psi,~\psi=\left\langle y,\mathbf{e}\left(
\varphi\right)  \right\rangle $ where $y_{1}=x_{1}/x_{0},y_{2}=x_{2}/x_{0}$
defined in $X\times\Sigma.$ By Lemma \ref{T} we have $N\left(  x,x^{\prime
}\right)  =0$ for any $x\neq x^{\prime}\in X$. Let \textrm{g }be the standard
spherical metric in $X.$ We have
\begin{align*}
\left\vert \nabla_{\mathrm{g}}\psi\right\vert ^{2}  &  =x_{0}^{-2}\left(
1+\left\langle y,\mathbf{e}\right\rangle ^{2}\right)  =\frac{1+\lambda^{2}%
}{x_{0}^{2}},\ \\
\int\frac{\mathrm{d}\varphi}{\left\vert \nabla_{\mathrm{g}}\psi\right\vert
^{2}}  &  =x_{0}^{2}\int_{0}^{2\pi}\frac{\mathrm{d}\varphi}{1+y^{2}\cos
^{2}\varphi}=\frac{2\pi x_{0}^{2}}{\sqrt{1+y^{2}}}=2\pi x_{0}^{3}%
\end{align*}
By Theorem \ref{R} for any function $f\in L_{2}\left(  X\right)  $%
\begin{equation}
f\left(  x\right)  =-\frac{1}{4\pi^{2}}\frac{1}{x_{0}^{3}}\int\int
\frac{Mf\left(  \lambda,\varphi\right)  \mathrm{d}\lambda\mathrm{d}\varphi
}{\left(  \lambda+\left\langle y,\mathbf{e}\left(  \varphi\right)
\right\rangle \right)  ^{2}} \label{17}%
\end{equation}
Choose coordinates $\theta,-\pi/2\leq\theta\leq\pi/2,\varphi$ in the dual
sphere $\Sigma$ so that $\sigma=\left(  \sin\theta,\cos\theta\cos\varphi
,\cos\theta\sin\varphi\right)  \in\Sigma.\ $We have$\ \lambda=\tan
\theta,\ \lambda+\left\langle y,\mathbf{e}\left(  \varphi\right)
\right\rangle =\left\langle \sigma,x\right\rangle /x_{0}\cos\theta$ and
$\left\vert \nabla_{\mathrm{g}}\psi\right\vert ^{-1}=x_{0}\cos\theta.$%
\[
Mf\left(  \lambda,\varphi\right)  =\int_{F\left(  \lambda,\varphi\right)
}\frac{f\mathrm{d}_{\mathrm{g}}s}{\left\vert \nabla_{\mathrm{g}}%
\psi\right\vert }=\cos\theta\int_{F\left(  \lambda,\varphi\right)  }%
x_{0}f\mathrm{d}_{\mathrm{g}}s=\cos\theta~R\left(  x_{0}f\right)  \left(
\sigma\right)
\]
We apply (\ref{17}) to the function $h\left(  x\right)  =x_{0}f\left(
x\right)  $ and rearrange this formula as follows%
\begin{equation}
f\left(  x\right)  =-\frac{1}{4\pi^{2}}\int_{0}^{2\pi}\int_{-\pi/2}^{\pi
/2}\frac{R_{\mathrm{g}}f\left(  \sigma\right)  \cos\theta\mathrm{d}%
\theta\mathrm{d}\varphi}{\left\langle \sigma,x\right\rangle ^{2}}=-\frac
{1}{2\pi^{2}}\int_{S_{+}}\frac{R_{\mathrm{g}}f\left(  \sigma\right)
\mathrm{d}\sigma}{\left\langle \sigma,x\right\rangle ^{2}} \label{12}%
\end{equation}
where $\mathrm{d}\sigma=\cos\theta\mathrm{d}\theta\mathrm{d}\varphi$ is the
area form in the hemisphere $S_{+}$. This formula coincides with original
Funk's result \cite{Fu} after a partial integration.

\section{Geodesic transform in Lobachevski plane}

Take a generating$\ $function $\Phi\left(  x;\lambda,\varphi\right)
=\lambda+\psi,\ \psi=-2\left(  \left\vert x\right\vert ^{2}+1\right)
^{-1}\left\langle x,\mathbf{e}\left(  \varphi\right)  \right\rangle
,\ -1<\lambda<1\ $in the unit disc $X=D$. The curves $F\left(  \lambda
,\varphi\right)  $ are geodesics in the Poincar\'{e} model of the Lobachevski
plane and $\mathrm{d}_{\mathrm{g}}s=2\left(  1-\left\vert x\right\vert
^{2}\right)  ^{-1}\mathrm{d}s$ is the hyperbolic metric. We have%

\[
\left\vert \nabla\psi\right\vert ^{2}=\left(  \frac{2}{1+\left\vert
x\right\vert ^{2}}\right)  ^{2}\left(  1-\psi^{2}\right)
\]
and%
\[
Mf\left(  \lambda,\varphi\right)  =\left(  1-\lambda^{2}\right)  ^{-1/2}%
\int\frac{1+\left\vert x\right\vert ^{2}}{2}f\left(  x\right)  \mathrm{d}%
s=\left(  1-\lambda^{2}\right)  ^{-1/2}\int\frac{1-\left\vert x\right\vert
^{4}}{4}f\left(  x\right)  \mathrm{d}_{\mathrm{g}}s,\
\]
where $\mathrm{d}s,\mathrm{d}_{\mathrm{g}}s$ means the Euclidean and
hyperbolic line element respectively. Further%
\[
D\left(  x\right)  =\frac{1}{2\pi}\int_{0}^{2\pi}\frac{\mathrm{d}\varphi
}{\left\vert \nabla\psi\left(  x,\varphi\right)  \right\vert ^{2}}=\frac
{1}{8\pi}\left(  1+\left\vert x\right\vert ^{2}\right)  ^{2}\int_{0}^{2\pi
}\frac{\mathrm{d}\varphi}{1-y^{2}\cos^{2}\varphi}=\frac{\left(  1+\left\vert
x\right\vert ^{2}\right)  ^{3}}{4\left(  1-\left\vert x\right\vert
^{2}\right)  }%
\]
since the integral in the right-hand side is equal to $2\pi\left(
1-\left\vert y\right\vert ^{2}\right)  ^{-1/2},\ y=2\left\vert x\right\vert
\left(  1-\left\vert x\right\vert ^{2}\right)  ^{-1}$ (e.g. Lemma \ref{I}). By
Theorem \ref{R} this yields%
\begin{equation}
f\left(  x\right)  =-\frac{1}{\pi^{2}}\frac{1-\left\vert x\right\vert ^{2}%
}{\left(  1+\left\vert x\right\vert ^{2}\right)  ^{3}}\int\int\frac{Mf\left(
\lambda,\varphi\right)  \mathrm{d}\lambda\mathrm{d}\varphi}{\Phi^{2}\left(
x;\lambda,\varphi\right)  } \label{4}%
\end{equation}
We apply this equation to the function $f^{\ast}\left(  x\right)  =4\left(
1-\left\vert x\right\vert ^{4}\right)  ^{-1}f\left(  x\right)  $ and get%
\[
f\left(  x\right)  =-\frac{1}{4\pi^{2}}\left(  1-\left\vert x\right\vert
^{2}\right)  ^{2}\int_{0}^{2\pi}\int_{-1}^{1}\frac{Rf\left(  \lambda
,\varphi\right)  \left(  1-\lambda^{2}\right)  ^{-1/2}\mathrm{d}%
\lambda\mathrm{d}\varphi}{\left(  \left(  1+\left\vert x\right\vert
^{2}\right)  \lambda-2\left\langle x,\mathbf{e}\left(  \varphi\right)
\right\rangle \right)  ^{2}}%
\]
where $R$ $f\left(  \lambda,\varphi\right)  =\int f\mathrm{d}_{\mathrm{g}}s$
and the integral is taken over a geodesic $F\left(  \lambda,\varphi\right)  $
in Lobachevski plane. This formula is equivalent to Helgason's theorem
\cite{He}.

Applying a map%
\[
x=\frac{y^{\prime}}{1+y_{0}},\ y^{\prime}=\left(  y_{1},y_{2}\right)
\]
we transform the Poincar\'{e} model to the Lorentz model in the "upper" sheet
hyperboloid
\[
Q_{+}=\{\left(  y_{0},y\right)  \in E^{3},y_{0}^{2}=1+y_{1}^{2}+y_{2}%
^{2},y_{0}>0\}
\]
with the metric induced from the Euclidean metric in $E^{3}.$ For any central
plane $P$ in $E^{3}$ the hyperbola $\gamma=P\cap Q_{+}$ is a geodesic curve in
$Q_{+}.$ Vice versa, the image in $Q_{+}$ of an arbitrary geodesic circle
$F\left(  \lambda,\varphi\right)  =\{\lambda=2\left(  1+\left\vert
x\right\vert ^{2}\right)  ^{-1}\left\langle x,\mathbf{e}\left(  \varphi
\right)  \right\rangle \}$ is contained in the plane $P=\{y;\lambda
y_{0}-\left\langle \mathbf{e}\left(  \varphi\right)  \mathbf{,}y^{\prime
}\right\rangle =0\}.$ This plane is orthogonal to a vector $\xi\doteq\left(
1-\lambda^{2}\right)  ^{-1/2}\left(  \lambda,-\mathbf{e}\left(  \varphi
\right)  \right)  $ which is contained in the one-sheet hyperboloid
$Q_{-}=\{\xi\in E^{3},\xi_{0}^{2}-\xi_{1}^{2}-\xi_{2}^{2}=-1\}.$ Therefore we
use notation $\gamma\left(  \xi\right)  =F\left(  \lambda,\varphi\right)
=P\cap Q_{+}$ for a geodesic in $Q_{+}$ and the hyperboloid $Q_{-}$
parametrizes all the geodesics in the Lorentz model. We can write$\ \Phi
\left(  x;\lambda,\varphi\right)  =\lambda-y_{0}^{-1}\left\langle
y,\mathbf{e}\left(  \varphi\right)  \right\rangle =y_{0}^{-1}\left(
1-\lambda^{2}\right)  ^{1/2}\left\langle \xi,y\right\rangle $ and have
$1-\left\vert x\right\vert ^{2}=2\left(  1+y_{0}\right)  ^{-1},\ 1+\left\vert
x\right\vert ^{2}=2y_{0}\left(  1+y_{0}\right)  ^{-1}.$ The equation (\ref{4})
is now read as follows%
\[
f\left(  y\right)  =-\frac{1}{4\pi^{2}}\int\int\frac{Rf\left(  \gamma\left(
\xi\right)  \right)  \left(  1-\lambda^{2}\right)  ^{-3/2}\mathrm{d}%
\lambda\mathrm{d}\varphi}{\left\langle \xi,y\right\rangle ^{2}}%
\]
Expressing the projective area form
\[
\omega\left(  \xi\right)  =\xi_{0}\mathrm{d}\xi_{1}\mathrm{d}\xi_{2}+\xi
_{1}\mathrm{d}\xi_{2}\mathrm{d}\xi_{0}+\xi_{2}\mathrm{d}\xi_{0}\mathrm{d}%
\xi_{1}%
\]
in terms of coordinates $\lambda$ and $\varphi$ we get $\omega\left(
\xi\right)  =\left(  1-\lambda^{2}\right)  ^{-3/2}\mathrm{d}\lambda
\mathrm{d}\varphi$. Finally
\[
f\left(  y\right)  =-\frac{1}{4\pi^{2}}\int_{Q_{+}}\frac{Rf\left(
\gamma\left(  \xi\right)  \right)  \omega\left(  \xi\right)  }{\left\langle
\xi,y\right\rangle ^{2}}%
\]
which makes complete similarity with\ (\ref{12}).

\textbf{Remark.} A formula of this form appears in \cite{GGG}, p.96 however
with misprinted numerical coefficient.

\section{Equidistant curves in Lobachevski plane}

Let $X=D$ be the open unit disc and $F=\{F\left(  \lambda,\varphi\right)  \}$
be the family of arcs in $D$ connecting two opposite points $\mathbf{e}\left(
\varphi-\pi/2\right)  $ and $\mathbf{e}\left(  \varphi+\pi/2\right)  $ on
$\partial D$ where $\lambda=\tan\omega,$ $\omega$ is the angular measure of
$F\left(  \lambda,\varphi\right)  .$ For any $\varphi$ the arcs $F\left(
\lambda,\varphi\right)  $ and $F\left(  \lambda^{\prime},\varphi\right)  $ are
equidistant with distance $2\left\vert \lambda-\lambda^{\prime}\right\vert .$
This family has generating function
\[
\Phi\left(  x;\lambda,\varphi\right)  =\lambda+\psi,\ \psi\left(
x,\varphi\right)  =-\frac{2\left\langle x,\mathbf{e}\left(  \varphi\right)
\right\rangle }{1-\left\vert x\right\vert ^{2}}%
\]
Check that this function fulfils the conditions of Theorem \ref{R}. A proof of
regularity is a routine. Further we have$\ \psi\left(  x,\varphi\right)
-\psi\left(  y,\varphi\right)  =\left\langle y^{\prime}-x^{\prime}%
,\mathbf{e}\left(  \varphi\right)  \right\rangle $ where $\ x^{\prime}=\left(
1-\left\vert x\right\vert ^{2}\right)  ^{-1}x,\ y^{\prime}=\left(
1-\left\vert y\right\vert ^{2}\right)  ^{-1}y.$ This is a first order
trigonometric polynomial with two real roots which implies $N\left(
x,y\right)  =0$ as $x\neq y.$ We have
\[
\left\vert \nabla\psi\left(  x,\varphi\right)  \right\vert ^{2}=\frac
{4}{\left(  1-\left\vert x\right\vert ^{2}\right)  ^{2}}\left(  1+z^{2}%
\cos^{2}\varphi\right)
\]
where $z=2\left\vert x\right\vert \left(  1-\left\vert x\right\vert
^{2}\right)  ^{-1}$ and
\[
D\left(  x\right)  =\frac{\left(  1-\left\vert x\right\vert ^{2}\right)  ^{2}%
}{8\pi}\int_{0}^{2\pi}\frac{\mathrm{d}\varphi}{1+z^{2}\cos^{2}\varphi}%
=\frac{1}{4}\frac{\left(  1-\left\vert x\right\vert ^{2}\right)  ^{3}%
}{1+\left\vert x\right\vert ^{2}}%
\]

\begin{corollary}
For any function $f$ with compact support in the unit disc a reconstruction
from $M_{\Phi}f$ is given by%
\[
f\left(  x\right)  =-\frac{1}{\pi^{2}}\frac{1+\left\vert x\right\vert ^{2}%
}{1-\left\vert x\right\vert ^{2}}\int_{0}^{2\pi}\int_{-1}^{1}\frac{M_{\Phi
}f\left(  \lambda,\varphi\right)  }{\left(  \left(  1-\left\vert x\right\vert
^{2}\right)  \lambda-2\left\langle x,\mathbf{e}\left(  \varphi\right)
\right\rangle \right)  ^{2}}\mathrm{d}\lambda\mathrm{d}\varphi
\]

\end{corollary}

Because of $\left\vert \nabla\Phi\left(  x;\lambda,\varphi\right)  \right\vert
=2\left(  1-\left\vert x\right\vert ^{2}\right)  ^{-1}\cos^{-1}\omega$ we can
express the Funk-Radon operator $M$ in terms of Euclidean arc integrals%
\[
Mf\left(  \lambda,\varphi\right)  =\frac{\cos\omega}{2}Rf^{\ast}\left(
\lambda,\varphi\right)  ,\lambda=\tan\omega,\ f^{\ast}\left(  x\right)
=\left(  1-\left\vert x\right\vert ^{2}\right)  f\left(  x\right)
\]
This yields%
\[
f\left(  x\right)  =-\frac{1+\left\vert x\right\vert ^{2}}{2\pi^{2}}\int
_{0}^{2\pi}\int_{-\pi/2}^{\pi/2}\frac{Rf\left(  \omega,\varphi\right)
\cos\omega\mathrm{d}\omega\mathrm{d}\varphi}{\left(  \left(  1-\left\vert
x\right\vert ^{2}\right)  \sin\omega-2\left\langle x,\mathbf{e}\left(
\varphi\right)  \right\rangle \cos\omega\right)  ^{2}}%
\]
This formula was stated in \cite{P3} by a different method.

\section{Photoacoustic geometries}

\textbf{Elliptical source curve. }D.\ Finch and Rakesh \cite{F} gave a formula
of type (\ref{3}) for reconstruction for the family of spheres centered on a
sphere. Another reconstruction formula was proposed by L. Kunyansky \cite{K};
it looks different however can be reduced to a form close to (\ref{3}) after
one-fold integration. We show that a reconstruction can be done for ellipse
and more general source curves by means of Theorem \ref{R}. F. Natterer
applied quite different method for reconstruction from the family of spheres
centered in an ellipsoid \cite{N}.

Take first a generating function $\Phi=\lambda+\psi$ where%

\[
\psi\left(  x,\varphi\right)  =\left\vert x-\mathbf{e}\left(  \varphi\right)
\right\vert ^{2},\mathbf{e}\left(  \varphi\right)  =\left(  \cos\varphi
e_{1},\sin\varphi e_{2}\right)  ,
\]
where $e_{1},e_{2}$ are arbitrary positive numbers. The source curve
$s=\mathbf{e}\left(  \varphi\right)  ,\varphi\in S^{1}$ is an ellipse with
half-axes $e_{1}$ and $e_{2}.$

\begin{proposition}
For any function $f$ with support in the ellipse $E=\{x;\left(  x_{1}%
/e_{1}\right)  ^{2}+\left(  x_{2}/e_{2}\right)  ^{2}\leq1\}$ a reconstruction
is given by%
\begin{align*}
f\left(  x\right)   &  =-\frac{\left\vert e\right\vert ^{2}-\left\vert
x\right\vert ^{2}}{\pi^{2}}\int\int\frac{Mf\left(  \lambda,\varphi\right)
}{\left(  \lambda-\left\vert x-\mathbf{e}\left(  \varphi\right)  \right\vert
^{2}\right)  ^{2}}\mathrm{d}\lambda\mathrm{d}\varphi\\
&  =-\frac{\left\vert e\right\vert ^{2}-\left\vert x\right\vert ^{2}}{\pi^{2}%
}\int\int\frac{Rf\left(  r,\varphi\right)  }{\left(  r^{2}-\left\vert
x-\mathbf{e}\left(  \varphi\right)  \right\vert ^{2}\right)  ^{2}}%
\mathrm{d}r\mathrm{d}\varphi
\end{align*}
where $e=\left(  e_{1},e_{2}\right)  $ and $\lambda=r^{2}$ since $Mf\left(
\lambda,\varphi\right)  =R$ $f\left(  r,\varphi\right)  /2r,$ $R$ $f\left(
r,\varphi\right)  =\int_{F\left(  r,\varphi\right)  }f\mathrm{d}s.$
\end{proposition}

\textbf{Proof. }We have%
\begin{align}
\psi\left(  x,\varphi\right)  -\psi\left(  y,\varphi\right)   &  =\left\vert
x-\mathbf{e}\left(  \varphi\right)  \right\vert ^{2}-\left\vert y-\mathbf{e}%
\left(  \varphi\right)  \right\vert ^{2}=2\left\langle y-x,\mathbf{e}\left(
\varphi\right)  \right\rangle +\left\vert x\right\vert ^{2}-\left\vert
y\right\vert ^{2}\nonumber\\
&  =2\left\Vert y-x\right\Vert _{e}\cos\left(  \varphi-\theta\right)
+\left\vert x\right\vert ^{2}-\left\vert y\right\vert ^{2} \label{13}%
\end{align}
where $\theta=\mathrm{\arg}\left(  y-x\right)  $ and $\left\Vert z\right\Vert
_{e}\doteq\left(  \left(  z_{1}e_{1}\right)  ^{2}+\left(  z_{2}e_{2}\right)
^{2}\right)  ^{1/2}.$ This norm and the dual norm $\left\Vert z\right\Vert
_{e}^{\ast}=\left(  \left(  z_{1}/e_{1}\right)  ^{2}+\left(  z_{2}%
/e_{2}\right)  ^{2}\right)  ^{1/2}$ fulfil the inequality
\begin{align*}
\left\vert \left\vert x\right\vert ^{2}-\left\vert y\right\vert ^{2}%
\right\vert  &  \leq\left\vert y_{1}-x_{1}\right\vert \left\vert y_{1}%
+x_{1}\right\vert +\left\vert y_{2}-x_{2}\right\vert \left\vert y_{2}%
+x_{2}\right\vert \\
&  \leq\left\Vert y-x\right\Vert _{e}\left\Vert y+x\right\Vert _{e}^{\ast}%
\leq2\left\Vert y-x\right\Vert _{e}%
\end{align*}
since$\ \left\Vert y+x\right\Vert _{e}^{\ast}<2$ for $x,y$ belonging to the
support of $f$. This implies that the trigonometric polynomial (\ref{13}) has
only simple real zeros. By Lemma \ref{T} $N\left(  x,y\right)  =0$ for
arbitrary $x\neq y\in E$. We have $\left\vert \nabla\psi\right\vert
=2\left\vert x-\mathbf{e}\left(  \varphi\right)  \right\vert $ and%
\[
D\left(  x\right)  =\frac{1}{8\pi}\int\frac{\mathrm{d}\varphi}{\left\vert
x\right\vert ^{2}+\left\vert e\right\vert ^{2}-2\left\vert x\right\vert
\left\vert e\right\vert \cos\left(  \varphi-\theta\right)  }=\frac{1}{4}%
\frac{1}{\left\vert e\right\vert ^{2}-\left\vert x\right\vert ^{2}}%
\]
and our statement follows. $\blacktriangleright$

A similar reconstruction from spherical means in $\mathbb{R}^{n}$ can be done
by means of Theorem \ref{N}.

\section{Isofocal hyperbolas and parabolas}

\textbf{Hyperbolas. }The equation $\lambda=\left\vert x\right\vert
-\varepsilon x_{1},\varepsilon>1$ defines a fold of the hyperbola
\[
\left(  \alpha x_{1}+\frac{\lambda}{\alpha}\right)  ^{2}-x_{2}^{2}%
=\frac{\lambda^{2}}{\alpha^{2}},\ \alpha=\sqrt{\varepsilon^{2}-1}%
\]
with a focus at the origin. The function $\Phi\left(  x;\lambda,\varphi
\right)  =\lambda+\psi\left(  x,\varphi\right)  ,$ $\psi\left(  x,\varphi
\right)  =\varepsilon\left\langle x,\mathbf{e}\left(  \varphi\right)
\right\rangle -\left\vert x\right\vert $ generates the family of all one-fold
hyperbolae with focuses in the origin. The first order trigonometric
polynomial%
\[
\psi\left(  x,\varphi\right)  -\psi\left(  y,\varphi\right)  =\varepsilon
\left\langle x-y,\mathbf{e}\left(  \varphi\right)  \right\rangle -\left\vert
x\right\vert +\left\vert y\right\vert
\]
has two real zeros if $x\neq y.$ since $\left\vert \left\vert x\right\vert
-\left\vert y\right\vert \right\vert <\varepsilon\left\vert x-y\right\vert . $
By Lemma \ref{T} the nucleus $N\left(  x,y\right)  $ vanishes for $y\neq x$.
We have $\left\vert \nabla\psi\right\vert ^{2}=1+\varepsilon^{2}%
-2\varepsilon\left\vert x\right\vert ^{-1}\left\langle x,\mathbf{e}%
\right\rangle $ and%
\[
D\left(  x\right)  =\frac{1}{2\pi}\int_{0}^{2\pi}\frac{\mathrm{d}\varphi
}{1+\varepsilon^{2}-2\varepsilon\left\vert x\right\vert ^{-1}\left\langle
x,\mathbf{e}\left(  \varphi\right)  \right\rangle }=\frac{1}{\varepsilon
^{2}-1}%
\]
(Lemma \ref{I}).

\begin{corollary}
For any smooth function $f$ with compact support in $\mathbb{R}^{2}$ the
equation holds%
\[
f\left(  x\right)  =-\frac{\varepsilon^{2}-1}{4\pi^{2}}\int_{0}^{2\pi}%
\int_{-\infty}^{\infty}\frac{Mf\left(  \lambda,\varphi\right)  \mathrm{d}%
\lambda\mathrm{d}\varphi}{\left(  \left\vert x\right\vert -\varepsilon
\left\langle x,\mathbf{e}\left(  \varphi\right)  \right\rangle \right)  ^{2}}%
\]

\end{corollary}

\textbf{Parabolas. }A parabola with a focus at the origin can be given by
$2px_{1}+x_{2}^{2}-p^{2}=0$ where $p\ $is a positive parameter. Set
$\lambda=p^{1/2}$ and write this equation in the form $\left(  x_{1}%
-\lambda^{2}\right)  ^{2}=\left\vert x\right\vert ^{2}$ which is equivalent to
$\left\vert x\right\vert +x_{1}=\lambda^{2}$ since $x_{1}\leq\lambda^{2}.$
Take a function $\Phi\left(  x;\lambda,\varphi\right)  =\lambda+\psi\left(
x,\varphi\right)  $ in $X=\mathbb{R}^{2}\backslash\{0\}$ where%
\[
\psi\left(  x,\varphi\right)  =-\sqrt{\left\vert x\right\vert +\left\langle
x,\mathbf{e}\left(  \varphi\right)  \right\rangle }=-\left(  2\left\vert
x\right\vert \right)  ^{1/2}\cos\left(  \varphi-\theta\right)  /2,\ \theta
=\mathrm{\arg~}x
\]
This function generates all parabolas with focus at the origin. For any $x\neq
y\in\mathbb{R}^{2}$ the trigonometric polynomial $\psi\left(  x,\varphi
\right)  -\psi\left(  y,\varphi\right)  $ is of order 1 and zero constant
term. Therefore it has two real roots and by Lemma \ref{T} the function
$N\left(  x,y\right)  $ vanishes if $y\neq x.$ We have
\[
\left\vert \nabla\psi\right\vert ^{2}=\frac{1}{2\left\vert x\right\vert
},\ Mf\left(  \lambda,\varphi\right)  =\int_{F\left(  \lambda,\varphi\right)
}\left(  2\left\vert x\right\vert \right)  ^{1/2}f\left(  x\right)
\mathrm{d}s,\ D\left(  x\right)  =2\left\vert x\right\vert
\]

\begin{corollary}
For any $L_{2}$-function $f$ with compact support in $\mathbb{R}^{2}%
\backslash\{0\}$ a reconstruction is given by%
\begin{align*}
f\left(  x\right)   &  =-\frac{1}{8\pi^{2}\left\vert x\right\vert }\int
_{0}^{2\pi}\int_{0}^{\infty}\frac{Mf\left(  \lambda,\varphi\right)
\mathrm{d}\lambda\mathrm{d}\varphi}{\left(  \lambda-\sqrt{\left\vert
x\right\vert +\left\langle x,\mathbf{e}\left(  \varphi\right)  \right\rangle
}\right)  ^{2}}\\
&  =-\frac{1}{4\pi^{2}\left(  2\left\vert x\right\vert \right)  ^{1/2}}%
\int_{0}^{2\pi}\int_{0}^{\infty}\frac{Rf\left(  \lambda,\varphi\right)
\mathrm{d}\lambda\mathrm{d}\varphi}{\left(  \lambda-\sqrt{\left\vert
x\right\vert +\left\langle x,\mathbf{e}\left(  \varphi\right)  \right\rangle
}\right)  ^{2}}%
\end{align*}

\end{corollary}

In \cite{D} Fourier coefficients of $f$ are reconstructed from $Mf$ by
Cormack's method.

\section{Cormack's curves}

For $\alpha>0$ an $\alpha$-curve is given in a plane by a equation%
\[
\lambda=r^{\alpha}\cos\left(  \alpha\left(  \theta-\varphi\right)  \right)
,\ \left\vert \theta-\varphi\right\vert \leq\pi/2\alpha
\]
where $x=r\mathbf{e}\left(  \theta\right)  $ and $\lambda,\varphi$ are
parameters \cite{Cor}. Suppose that $\alpha=k$ is natural and take the family
of curves generated by a function $\Phi=\lambda+\psi\left(  x,\varphi\right)
,$ $\psi\left(  x,\varphi\right)  =-r^{\alpha}\cos\left(  \alpha\theta
-\varphi\right)  .$ Note that in the case $k=2$ the curve $\psi\left(
x,\varphi\right)  =-\lambda$ is a hyperbola. The function $\psi$ always admits
a symmetry group $\mathbb{Z}_{k}$ acting by rotation of the plane by the angle
$2\pi/k$ about the origin. Consider a quotient manifold $X=\mathbb{R}%
^{2}\backslash\{0\}/\mathbb{Z}_{k};$ lifting of $\Phi$ to $X$ is a regular
generating function. We have $\psi\left(  x,\varphi\right)  =\operatorname{Re}%
\left(  \mathrm{\exp}\left(  -i\varphi\right)  \left(  x_{1}+ix_{2}\right)
^{k}\right)  $ and
\[
\psi\left(  x,\varphi\right)  -\psi\left(  y,\varphi\right)
=\operatorname{Re}\mathrm{\exp}-i\varphi\left(  \left(  x_{1}+ix_{2}\right)
^{k}-\left(  y_{1}+iy_{2}\right)  ^{k}\right)
\]
is a first order trigonometric polynomial with two real zeros for arbitrary
$x\neq y.$ This implies $N\left(  x,y\right)  =0$. Further we have
\[
\left\vert \nabla\psi\left(  x,\varphi\right)  \right\vert ^{2}=k^{2}%
\left\vert x\right\vert ^{2k-2},\ D\left(  x\right)  =\frac{1}{k^{2}\left\vert
x\right\vert ^{2k-2}},\ Mf\left(  \lambda,\varphi\right)  =\frac{1}{k}%
\int_{\lambda=\psi\left(  x,\varphi\right)  }\frac{f\left(  x\right)
\mathrm{d}s}{\left\vert x\right\vert ^{k-1}}%
\]
and the reconstruction formula reads: for any $\mathbb{Z}_{k}$-invariant
function $f$ with compact support in $\mathbb{R}^{2}\backslash\{0\}:$%
\[
f\left(  x\right)  =-\frac{k\left\vert x\right\vert ^{k-1}}{4\pi^{2}}\int
_{0}^{2\pi}\int_{-\infty}^{\infty}\frac{\ f\left(  \lambda,\varphi\right)
\mathrm{d}\lambda\mathrm{d}\varphi}{\left(  \lambda-\operatorname{Re}%
\mathrm{\exp}\left(  -i\varphi\right)  \left(  x_{1}+ix_{2}\right)
^{k}\right)  ^{2}}%
\]
The family of $\beta$-curves $\lambda=r^{-\beta}\cos\beta\left(
\theta-\varphi\right)  $ with natural $\beta$ is studied in the same way.

\section{Higher dimensions}

Let $\Phi$ be a generating function defined in $X\times\Sigma$ where
$X\subset\mathbb{R}^{n}$ is a\ \ smooth manifold of dimension $n$ with a
Riemannian metric $\mathrm{g}$,$\ \Sigma=\mathbb{R}\times S^{n-1}$ and
$\Phi\left(  x;\lambda,\omega\right)  =\lambda+\psi\left(  x,\omega\right)  $,
$\lambda\in\mathbb{R}$, $\omega\in S^{n-1}.$ We say that a generating function
$\Phi$ is regular if a $n+1\times n+1$ determinant like (\ref{15}) does not
vanish in $F=\{\Phi=0\}$ and conjugated points are absent, see \cite{P2} for
details. The Funk-Radon transform is defined by%
\[
M_{\Phi}f\left(  \lambda,\varphi\right)  =\int_{F\left(  \lambda,\phi\right)
}\frac{f\mathrm{d}V}{\mathrm{d}_{x}\Phi}%
\]
where$\ \mathrm{d}V$ means the Riemannian volume form and $F\left(
\lambda,\varphi\right)  =\left\{  x:\Phi\left(  x;\lambda,\varphi\right)
=0\right\}  .$ The result of Sec.\ 3 is generalized as follows. Define a
function%
\[
N\left(  x,y\right)  =\operatorname{Re}\int_{S^{n-1}}\frac{\mathrm{d}\omega
}{\left(  \psi\left(  x,\omega\right)  -\psi\left(  y,\omega\right)  \pm
i0\right)  ^{n}},x\neq y
\]
where $\mathrm{d}\omega$ the Euclidean volume form in $S^{n-1}$ and%
\[
D\left(  x\right)  =\frac{1}{\left\vert S^{n-1}\right\vert }\int_{S^{n-1}%
}\frac{\mathrm{d}\omega}{\left\vert \nabla_{\mathrm{g}}\psi\left(
x,\omega\right)  \right\vert ^{n}}%
\]

\begin{theorem}
\label{N}If a regular generating function $\Phi$ is real analytic in $X$ and
$N\left(  x,y\right)  =0$ for any $x\neq y\in X$ then a reconstruction of type
(\ref{3}) exists for an arbitrary function $f\in L_{2}\left(  X\right)  $ with
compact support. For even $n$ it reads%
\[
f\left(  x\right)  =\frac{1}{\left(  2\pi i\right)  ^{n}D\left(  x\right)
}\operatorname{Re}\int_{\Sigma}\frac{M_{\Phi}f\left(  \lambda,\omega\right)
\mathrm{d}\lambda\mathrm{d}\omega}{\left(  \Phi\left(  x;\lambda
,\omega\right)  +i0\right)  ^{n}}%
\]
For odd $n$ it is%
\[
f\left(  x\right)  =\frac{1}{2\left(  2\pi i\right)  ^{n-1}D\left(  x\right)
}\int_{\Sigma}\delta^{\left(  n-1\right)  }\left(  \Phi\left(  x;\lambda
,\omega\right)  \right)  M_{\Phi}f\left(  \lambda,\omega\right)
\mathrm{d}\lambda\mathrm{d}\omega
\]
The integrals converge in $L_{2}\left(  X\right)  _{\mathrm{loc}}.$
\end{theorem}

A proof will be given elsewhere.

\end{document}